\documentclass[12pt,twoside]{article}

\usepackage{amsgen,amsmath,amstext,amsbsy,amsopn,amsfonts,amssymb}

\usepackage{jltmac2ePrep}

\renewcommand{\version}{{\rm Version of~}\date} 
                         
\renewcommand{\title}{\centerline{Uniqueness in
the Kashiwara-Vergne conjecture}}  
\renewcommand{\author}{A. Alekseev and E. Petracci}                 
\renewcommand{\lastname}{Alekseev and Petracci}                  
\renewcommand{\entries}{
\nextref{AlM}{Alekseev, A., and  E. Meinrenken} {\it 
	Poisson	geometry and the
        Kashiwara-Vergne conjecture}{ C. R. Acad. Sci., Paris,
        Ser. I   {\bf 335} (2002), 723--728}

\nextref{AlM3}{Alekseev, A., and E. Meinrenken}{\it On the 
Kashiwara-Vergne 
conjecture}
        {Preprint {\bf  math.QA/0506499}}

\nextref{Bou}{Bourbaki, N.}{``El\'ements de math\'ematiques, 
groupes et alg\`ebres de Lie, chapitres 2 
et 3''}{Hermann,  1972}

\nextref{KaV}{Kashiwara, M., and  M. Vergne}{\it The 
	Campbell-Hausdorff
        Formula and Invariant Hyperfunctions}{Inventiones math. 
        {\bf 47} (1978), 249-272}

\nextref{Tor}{Torossian, C.}  {\it Sur la conjecture combinatoire
         de Kashiwara-Vergne}{  J. Lie Theory {\bf 12} (2002), no. 
	2, 597--616}      

\nextref{Ver}{Vergne, M.} {\it Le centre de l'alg{\`e}bre enveloppante 
et        la formula de Campbell-Hausdorff}{ C. R. Acad. Sci., Paris,
        Ser. I Math.  {\bf 329} (1999), 767--772}
}

\renewcommand{\address}{Anton Alekseev

University of Geneva
\\
Section of Mathematics
\\
C.P. 64, 2-4 rue du Li\`evre
\\
CH-1211 Gen\`eve 4
\\
\makebox[1in][l]{Anton.Alekseev{@}math.unige.ch} 
}

\renewcommand{\addresstwo}{Emanuela Petracci           

University of Geneva
\\
Section of Mathematics
\\
C.P. 64, 2-4 rue du Li\`evre
\\
CH-1211 Gen\`eve 4
\\
\makebox[1in][l]{Emanuela.Petracci{@}math.unige.ch}

}

\DeclareMathAlphabet\EuFrak{U}{euf}{m}{n}
\SetMathAlphabet\EuFrak{bold}{U}{euf}{b}{n} 
\newcommand{\gG} 
        {{\EuFrak g} }
\newcommand {\tr}
	{{\rm tr}} 
\newcommand {\campo} 
        {\mathbb{K}}
\newcommand {\ad}
        { {\rm ad} }
\newcommand {\Reali}
        {\mathbb{R}}
\newcommand {\complessi}
        {\mathbb{C}}

\newcommand {\Naturali}
        {\mathbb{N}}  

\begin{document}
\firstpage

\begin{abstract} 
We prove that a universal symmetric solution of 
the Kashiwara-Vergne conjecture is unique up to order one. 
In the Appendix by the second author, this result is used
to show that solutions of the Kashiwara-Vergne 
conjecture for quadratic Lie algebras existing in the literature 
are not universal. 
\end{abstract} 

\section{Introduction}
Let $G$ be a Lie group over $\campo = \Reali$ (or 
$\campo=\complessi$), 
and $\gG$ its Lie algebra. There exists an open neighbourhood $\gG_0$ 
of the origin $0 \in \gG_0\subseteq\gG$ such that the restriction to 
$\gG_0$ of the exponential map $exp :\gG\rightarrow  G$ is 
an analytic diffeomorphism. We denote by $\ln:exp(\gG_0)\rightarrow \gG_0$ 
the inverse map and by 
$$
   \varphi_1(t) 
:= \frac{t}{e^t-1}
 = 1-\frac{1}{2}t + o(t)
$$  
the generating series of Bernoulli numbers. It is convenient 
to have a separate notation for the function
$\psi(t):= -(\varphi_1(t)-1)/2$.

In {\bf\cite{KaV}}, Kashiwara and Vergne  put forward 
the following conjecture on the properties of the 
Campbell-Hausdorff series:

\vskip0.5cm
\noindent 
{\bf Kashiwara-Vergne conjecture.}    
There exists a pair of $\gG$-valued analytic functions $A$ and $B$ 
defined 
on an open subset  $U \subset \gG\times\gG$  containing $(0,0)$, such that 
$A(0,0)=B(0,0)=0$, and  
for any $(X,Y)\in U$ one has 
{\small
\begin{eqnarray}
\label{primaE}
&&\ln(exp(Y)exp(X))
- X -Y
= (id - e^{-\ad X})A(X,Y)
+ (e^{\ad Y}-id)B(X,Y),
\quad \quad \quad \quad 
\\\nonumber\\
\nonumber
&& \tr\left(
	  \ad X\circ \delta_1A(X,Y)
	+ \ad Y\circ \delta_2B(X,Y)
  \right)=
\\
\label{secondaE}
&&
\quad \quad \quad 
\quad \quad \quad 
\quad \quad \quad 
= \tr\big( \psi(\ad X) 
	+ \psi(\ad Y)
	- \psi\left( \ad \ln\big(exp(X)exp(Y)\big)
	  \right)
    \big)
,
\end{eqnarray}
where $\delta_1A(X,Y), \delta_2B(X,Y) \in {\rm End}(\gG)$ are
defined as follows,
$$
\delta_1A(X,Y): Z \mapsto \frac{d}{dt} A(X+tZ,Y)|_{t=0} \, , \,
\delta_2B(X,Y): Z \mapsto \frac{d}{dt} A(X,Y+tZ)|_{t=0} \, .
$$
}
%
\vskip 0.3cm

Sometimes this statement is referred to as the
`combinatorial Kashiwara-Vergne conjecture' 
(see {\em e.g.} {\bf\cite{Tor}}). This conjecture
was established for solvable Lie algebras 
in {\bf\cite{KaV}} and for quadratic Lie algebras
in {\bf\cite{Ver}}. Recently, the general case
was settled in {\bf\cite{AlM3}} based on the 
earlier work {\bf \cite{Tor}}.

We denote by $\campo[[t]]$ and  $\campo[t]$ the ring of formal 
power series and  the ring of polynomials, respectively.
We call a solution of the Kashiwara-Vergne conjecture 
{\em universal} if $A$ and $B$ are given by series 
in Lie polynomials of the variables $X$ and $Y$:
\begin{eqnarray*}
A(X,Y) = \rho X + \beta(\ad X)(Y) + o(Y)
\\
B(X,Y) = \alpha X + \gamma(\ad X)(Y) + o(Y)
\end{eqnarray*}
with     $\beta(t)$, $\gamma(t)\in\campo[[t]]$,  
$\alpha, \rho\in\campo$,  
and  both $o(Y)$ are of type 
$$
      o(Y) 
  \in \sum_{k\geq 2}
      \sum_{
	\begin{array}{c}
	j_1,..., j_k\geq 0
	\\
	j_{k-1}<j_k
	\end{array}
	   }
\campo \, 
      \ad_{(\ad X)^{j_1}Y}
\circ \cdots
\circ \ad_{(\ad X)^{j_{k-1}}Y}
\circ (\ad X)^{j_k}(Y)
.
$$
If $(A,B)$ is a universal solution, the coefficients of the
Taylor expansions of $A$ and $B$ are the same for all Lie
algebras over $\campo$.

The set of solutions of the Kashiwara-Vergne conjecture carries 
a natural $\mathbb{Z}/2\mathbb{Z}$-action,
$$
(A(X,Y), B(X,Y)) \mapsto (B(-Y,-X), A(-Y,-X)) .
$$
A solution is called {\em symmetric} if it is stable with respect to
this action. Averaging of any solution produces a symmetric
solution. Hence, without loss of generality we can restrict
our attention to symmetric solutions. 

It is well-known (see {\em e.g.} {\bf\cite{Tor}}) that $\alpha, \rho$
and $\beta(t)$ are uniquely determined by the Kashiwara-Vergne
equations and by the symmetry condition. 
In this note we prove the uniqueness statement for the function $\gamma(t)$.
Thus, the symmetric universal solution of the Kashiwara-Vergne
conjecture is unique up to order one in $Y$.

In the Appendix by the second author,
this result is applied to show that solutions of the
Kashiwara-Vergne conjecture for quadratic Lie algebras
obtained in {\bf\cite{Ver}} and {\bf\cite{AlM}} are
not universal.

\section{Preliminaries}
\label{Sdue}

In this Section, we collect some elementary properties
of Lie algebras.

\begin{Remark} {\bf (Free Lie algebras with two 
generators).}
\label{bella}
We denote by $L_\campo(x,y)$ the free Lie 
$\campo$-algebra with generators $x$ and $y$. 
In this section  we use the Hall  basis $H$ of $L_\campo(x,y)$ 
defined in 
{\bf\cite{Bou}} 
(Definition 2,   page 27). 

$H$ consists of Lie words with the following order relation: 
$x,y \in H$ and $x<y$; if the number of Lie brackets in $a\in H$
is smaller than the number of Lie brackets in $b\in H$ then $a<b$;
and we omit the description of the order relation for $a$ and $b$
of equal length. The basis $H$ is built inductively starting
with $x, y, [x,y]$, and one adds the elements of the form
$[a,[b,c]]$ such that $a,b,c,[b,c]\in H$, $b\leq a\leq [b,c]$,
and $b<c$.

\end{Remark}

Using  the definition of $H$ we can prove by induction that 
$$
\forall n\geq 0\ \ \ \  (\ad x)^n (y)\in H.
$$
In fact, the cases $n=0$ and $n=1$ are 
trivial, and  for $n\geq 1$ we use $(\ad x)^{n+1}(y)= [x, [x,(\ad 
x)^{n-1}(y)]]$. Furthermore,
\begin{equation}
\label{tm21}
\forall n\geq 1
\ \ \ 
\left\{
[ (\ad x)^j(y)
, (\ad x)^{n-j}(y)
]
, 0\leq j<n-j, j\leq n-1
\right\}
\subset H
.
\end{equation} 
Here it is sufficient to observe that
$
  (\ad x)^{n-j}(y)
= [x,(\ad x)^{n-1-j}(y)]
$.

\begin{Proposition}
\label{lem21}
Let  $\xi(t)\in\campo[[t]]$. 
The following statements are equivalent:
\begin{enumerate}
\item[i)]
  for any Lie $\campo$-algebra $\gG$ we have 
$\xi(\ad X)(Y) = 0$ \ $\forall X,Y\in\gG$;  
\item[ii)] $\xi(t)=0$.
\end{enumerate}
\end{Proposition}
\begin{Proof}
It is sufficient to show that $i)$ implies $ii)$.  
Let $n\in\Naturali$. 
By rescaling $X\mapsto tX$ and applying
 $\frac{d^n}{dt^n}|_{t=0}$ 
we get
$\xi_n(\ad X)^n (Y)  =0$. Choosing $\gG=L_\campo(x,y)$, $X=x$ and $Y=y$ we 
get $\xi_n=0$.
\end{Proof}

The following will be a very useful notation.  

\begin{Definition} 
Let $W,X,Y\in\gG$.  
For any pair $i,j\in\Naturali$, we set 
$$
   (t^iu^j: [W,X])_Y
:= [(\ad Y)^i (W), (\ad Y)^j (X)] 
.
$$
This notation is extended by linearity 	
to any formal power series $\xi(t,u)\in\campo[[t,u]]$. Then 
$(\xi(t,u): [W,X])_Y\in \gG[[\gG]]$ is a formal power 
series with coefficients in $\gG$.
\end{Definition}

\begin{Remark}
\label{oss22}
\begin{enumerate}
\item[$i)$]  
$  
   (t^iu^j: [W,X])_Y
=- (u^it^j: [X,W])_Y
$.
\item[$ii)$] Jacobi's identity gives $(t+u: [W,X])_Y = 
(\ad Y)([W,X])$.
\end{enumerate} 
\end{Remark}

\begin{Proposition}
\label{lem23}
Let $\xi(u)=-\xi(-u)$ be a series in $\campo[[u]]$. The  following 
statements are equivalent:
\begin{enumerate}
\item[i)] for any Lie $\campo$-algebra $\gG$ we have 
$(\xi(u): [X,X])_Y = 0 \ \ \forall X,Y\in\gG$; 
\item[ii)] 
$\xi(u)=0$.
\end{enumerate}
\end{Proposition}
\begin{Proof}
Recall that  
$ 
      (\xi(u): [X,X])_Y 
    = [X, \xi(\ad Y)(X)]
$.
Similar to the proof of Proposition \ref{lem21}, it is 
sufficient to show that in the free Lie algebra $L_\campo(X,Y)$
we have $ [X,(\ad Y)^{2i+1}X] \neq 0$  for any $i\in\Naturali$
. Indeed, if we
rename $X=y, Y=x$, the elements $[y, (\ad x)^{2i+1}y]$
belong to the basis $H$ and, hence, are non-vanishing.
\end{Proof}

		
Every  formal power series 
$\xi(t,u)\in\campo[[t,u]]$ can be split  into the sum of its symmetric 
and
skew-symmetric parts:
$$
  \xi(t,u)
= \xi(t,u)_{skew} 
+ \xi(t,u)_{sym} 
= \frac{\xi(t,u)-\xi(u,t)}{2} 
+ \frac{\xi(t,u)+\xi(u,t)}{2}
.
$$

\begin{Proposition}
\label{cor21}
 Let $\xi(t,u)\in\campo[[t,u]]$. The following statements are 
equivalent:  
\begin{enumerate}
\item[i)]
 for any Lie $\campo$-algebra $\gG$ we have 
$(\xi(t,u): [X,X])_Y = 0$ $\forall X,Y\in\gG$; 
\item[ii)]
 $\xi(t,u)=\xi(u,t)$.
\end{enumerate}
\end{Proposition}
\begin{Proof}
By skew-symmetry of the Lie bracket, $(\xi(t,u)_{sym}: [X,X])_Y = 0$
for any $\xi$, and $(\xi(t,u): [X,X])_Y =(\xi(t,u)_{skew}
: [X,X])_Y$.

Let $\xi(t,u)$ be a formal power series with vanishing symmetric part.
Then, it can be written as
$$
\xi(t,u) = \frac{1}{2} \sum_{n=1}^\infty \sum_{0\leq j<n-j} 
\xi_{n,j} (t^j u^{n-j} - u^j t^{n-j}).
$$
Suppose that 
$$
(\xi(t,u): [X,X])_Y 
\equiv \sum_{n=1}^\infty \sum_{0\leq j<n-j}
\xi_{n,j} [(\ad Y)^j X, (\ad Y)^{n-j} X]  =0 
$$
for every Lie $\campo$-algebra $\gG$ and every $X,Y\in\gG$. 
By rescaling $X\mapsto t X$ and then applying the $n$-th 
derivative 
in $t$  we get 
$$
\forall n\geq 1
 \ \ \ \ 
  \sum_{0\leq j<n-j}
  \xi_{n,j} 
  [(\ad Y)^j X, (\ad Y)^{n-j} X]  
= 0 
.
$$ 
Then we choose 
$\gG=L_\campo(x,y)$, $X=y$ and $X=y$. 
Since all Lie words in the sum are linearly independent (recall 
property (\ref{tm21}))
this implies $\xi_{n,j}=0$ for all $n,j$ and $\xi(t,u)=0$.
 \end{Proof}

\begin{Lemma}
\label{l14}

In the Lie algebra $L_\campo(x,y)$ we have 
$$
\forall n\in\Naturali
,\ \ \ \
       (u^{2n+1}:[y,y])_x
\notin span_\campo
       \{ \big( (t+u)t^lu^m : [y,y]\big)_x
	  | l,m\in\Naturali
       \}
.
$$
\end{Lemma}
\begin{Proof}
We want to show that  
$ 
       (u^{2n+1}:[y,y])_x
\notin span_\campo
       \{ ((t+u)t^lu^{2n-l}: [y,y])_x
	| 0\leq l\leq 2n
       \}
.
$ 
If $n=0$ this statement is obvious. Let $n\geq 1$, and  
suppose that we can some find coefficients  $c_j\in\campo$ such 
that 
\begin{equation}
\label{ultima}
      (u^{2n+1}:[y,y])_x
    = \sum_{j=0}^{2n}
      c_j
      ((t+u)u^{2n-j}t^j:[y,y])_x
.
\end{equation}
Let 
$
   \xi(t,u) 
:= u^{2n+1}
 - \sum_{j=0}^{2n}
   c_j
   (t+u)
   u^{2n-j}
   t^j
$, then   
identity (\ref{ultima}) can be written as  $(\xi(t,u):[y,y])_x=0$. The 
universal 
property of a free Lie algebra allows to apply Proposition \ref{cor21}, 
so $\xi(t,u)=\xi(u,t)$. This  means that $u^{2n+1}-t^{2n+1} = 0$ 
modulo 
$(t+u)$, and this is a contradiction.
\end{Proof}

\begin{Remark}
\label{fine}
Here we explain that Propositions \ref{lem21}, 
\ref{lem23}, and  \ref{cor21} 
 still 
apply  if we 
restrict to finite-dimensional Lie algebras.  

In their proofs, at some point  we choose  $\gG$ equal to the  
free Lie algebra $L_\campo(x,y)$. 
Let  $N\geq 2$. We introduce  $\gG_N:= L_\campo(x,y)/I_N$, 
where 
$I_N$ is an ideal of $L_\campo(x,y)$ such that $\gG_N$ is 
an $N$-nilpotent Lie algebra. In particular $\gG_N$  is a 
finite-dimensional Lie algebra with basis $H/I_N$. To modify 
the  proofs  it is sufficient 
to replace 
$L_\campo(x,y)$ with   $\gG_N$, for a good choice of $N$: $n\leq N-1$ in 
Proposition \ref{lem21},  $2i+2\leq N-1$ in Proposition \ref{lem23}, 
and $n\leq N-2$ in Proposition \ref{cor21}.
\end{Remark}
\noindent
\newline
In the previous theorems we do not use Lie groups. We 
end this section by computing some derivatives of the exponential map of a 
Lie group $G$ with Lie algebra $\gG$. 

Let $g\in G$. We use  the notation $R_g :G\rightarrow G$ for 
the right translation. In the following lemma we denote by  $1_G$ the 
group unit of $G$.
 
\begin{Lemma}
\label{lemma21}
Let  $X,Y\in\gG$, then 
\begin{enumerate}
\item[i)] 
$
  d(\ln\circ R_{exp(X)})_{1_G}
= \varphi_1(\ad X)  
$,
\item[ii)]
$ 
  \frac{d}{ds}
  \ln(exp(sY)exp(X))
= \varphi_1(\ad _{\ln(exp(sY)exp(X)})
  (Y)
$,
\item[iii)] 
$
  \frac{d^2}{ds^2}|_{s=0}
  \ln(exp(sY)exp(X))
= \left(
     \frac{\varphi_1(t+u)-\varphi_1(u)}{t}
     \varphi_1(t)
   : [Y,Y]
  \right)_X
$.
\end{enumerate}
\end{Lemma}
\begin{Proof*}
$i)$ The formula of this differential is a consequence of the well-known 
formula for the differential of 
the exponential map: 
\begin{equation}
\label{difExp}
  d(exp)_X 
= (dL_{exp(X)})_{1_G}
\circ \varphi_1(-\ad X)^{-1}
.
\end{equation}
$ii)$  Using part $i)$ and $R_{exp(X)}R_{exp(-X)}=id$  we get 
\begin{eqnarray*}
&&      \frac{d}{ds}   
      \ln 
\circ R_{exp(X)}(exp(sY))
\equiv d(\ln\circ R_{exp(X)})
     _{exp(sY)}
      \frac{d}{ds}(exp(sY))=
\\
&&
    = \varphi_1(\ad _{\ln(exp(sY)exp(X))})
\circ (dR_{exp(-sY)})_{exp(sY)}
      \frac{d}{ds}(exp(sY))
.
\end{eqnarray*}
Formula (\ref{difExp}) gives 
$
  \frac{d}{ds}(exp(sY))
=  d(L_{exp(sY)})_{1_G}(Y)
$, so 
$$
  d(R_{exp(-sY)})_{exp(sY)}
  \frac{d}{ds}(exp(sY))
= e^{s\ \ad Y}(Y)
= Y
.
$$
$iii)$ Using $ii)$,   $i)$ and a direct calculation we get 
\begin{eqnarrayqed*}
&&      \frac{d^2}{ds^2}|_{s=0}   
      \ln 
\circ R_{exp(X)}(exp(sY))=
\\
&&\quad\quad
    = \left(
        \frac{\varphi_1(t+u) - \varphi_1(u)}
        {t}
      : [\frac{d}{ds}|_{s=0}\ln\circ R_{exp(X)}(exp(sY)),Y]
      \right)_X
\\
&&\quad\quad
\equiv  \left(
        \frac{\varphi_1(t+u) - \varphi_1(u)}
        {t}
        : [\varphi_1(\ad X)(Y),Y]
        \right)_X
.
\end{eqnarrayqed*}
\end{Proof*}

\section{The Campbell-Hausdorff series}
\label{Stre}
In this section we only make use of Equation (\ref{primaE}), and we derive
formulas for
$\beta(t)$, $\gamma(t)_{odd}$, and  
$ A_2(X,Y) 
:= \frac{d^2}{ds^2}|_{s=0}
  A(X,sY)
$.

\begin{Theorem}
\label{tm31bis}
A universal solution of the Kashiwara-Vergne conjecture has
$$
   \beta(t)
 = \beta_\alpha(t)
\equiv\varphi_1(-t)
   \left(
	\frac{\varphi_1(t)-1}{t} + \alpha
   \right)
.
$$
\end{Theorem}
\begin{Proof} 
In (\ref{primaE}) we rescale  $Y$ by $tY$ and we compute the 
derivative in $t=0$:
$$
  d(\ln\circ R_{exp(X)})_1(Y)
- Y
= (id - e^{-\ad X})
\circ  \beta(\ad X)(Y)
- \alpha (\ad X)(Y)
.
$$
Using Lemma \ref{lemma21} part $i$ and  
$
   (id - e^{-\ad X})
\equiv \varphi_1(-\ad X)^{-1}
\circ\ad X
$, we get 
 \begin{eqnarray*}  
&&
  \varphi_1(-\ad X)
\circ\big(\varphi_1(\ad X)
   - id
   + \alpha \ad X
  \big)
  (Y)
= \ad X
\circ  \beta(\ad X)
  (Y)
. 
\end{eqnarray*}
As this identity has to be verified for any $\campo$-Lie algebra $\gG$ and 
for any 
$X,Y\in\gG$,  
 Proposition \ref{lem21} implies   
$
  \varphi_1(-t)
  (  \varphi_1(t)
   - 1
   + \alpha t
  )
= t \beta(t)
$.
  \end{Proof}

The following theorem   uses the notation 
$$  A_2(X,Y)
= (\pi(t,u):[Y, Y])_X
,
$$
where the formal power series 
$\pi(t,u)\in\campo[[t,u]]$ is  skew-symmetric (i.e 
$\pi(t,u)= \pi(t,u)_{skew}$).
\begin{Theorem}
\label{tm31}
A universal solution of the Kashiwara-Vergne  conjecture  has  
\begin{enumerate}
\item[$i)$]
$
  \gamma(t)_{odd}
= \gamma_\alpha(t)_{odd}
:= \frac{\alpha}{2}
  t
+ \frac{1}{2}
  \left(
	\frac{\varphi_1(t)-1}{t}
	\varphi_1(-t)
  \right)_{odd}
,
$
\item[$ii)$]
$
  \pi(t,u)= \pi_{\gamma(t)}(t,u)
:=-\varphi_1(t+u)
  \frac{ \left(
	 \frac{\varphi_1(t+u)-\varphi_1(u)}{t}
         \varphi_1(t)
       + \alpha u
       - 2\gamma(u)
         \right)_{skew}
  }
  {t+u}
.
$
\end{enumerate}
\end{Theorem}
\begin{proof}
 By rescaling $Y \mapsto s Y$ and then applying the second
derivative in $s$ to Equation (\ref{primaE}) we obtain,
\begin{equation} \label{d^2}
  \frac{d^2}{ds^2}|_{s=0}
  \ln(exp(sY)exp(X))
= (id-e^{-\ad X})A_2(X,Y)
+ (-\alpha u + 2\gamma(u)
  : [Y, Y])_X .
\end{equation}
Here we have used that
\begin{eqnarray*}
&&\frac{d^2}{ds^2}|_{s=0}
  (e^{(\ad Y)s}-id)B(X,s Y)
= \alpha 
  (\ad Y)^2(X)
+ 2\ad Y\circ \gamma(\ad X)(Y)=
\\
&&
\quad\quad\quad\quad\quad\quad
\quad\quad\quad\quad
\quad\quad\quad
= (-\alpha u + 2\gamma(u): [Y, Y])_X
. 
\end{eqnarray*}

Let 
$s(t,u)
:= \frac{\varphi_1(t+u)-\varphi_1(u)}{t}
   \varphi_1(t)
 + \alpha u
 - 2\gamma(u)
$. 
The comparison of part $iii)$ of Lemma \ref{lemma21} 
with Equation (\ref{d^2}) gives
\begin{eqnarray*}
 \left(
     s(t,u)
   : [Y,Y]
  \right)_X
=  (1-e^{\ad X})A_2(X,Y) 
\equiv  \left( (1-e^{t+u})\pi(t,u):[Y,Y]
	\right)_X
\end{eqnarray*}
for any $X,Y\in\gG$ and any Lie $\campo$-algebra $\gG$. 

Let 
$
   g(t,u)
:= s(t,u) 
-  (1-e^{t+u})\pi(t,u)
$. 
 Proposition \ref{cor21} gives $g(t,u)_{skew}=0$, 
in particular 
$$
s(t,u)_{skew}=(1-e^{t+u})\pi(t,u)
.
$$ 
Putting $t+u=0$ we get  statement 
$i)$. To get statement $ii)$ it is sufficient to remark that 
$
(1-e^{t+u})= -\varphi_1(t+u)^{-1}(t+u)$.
\end{proof}

\section{The equation with traces}
\label{Squattro}
In this section we derive formulas for $\rho$ and $\gamma(t)$. 
We begin with a technical remark. 

\begin{Remark}
\label{os41}
Let $\lambda,\mu\in\campo\setminus\{0\}$ be two distinct numbers, and 
$\gG_{\lambda,\mu} =  \campo a \oplus \campo 
b\oplus \campo c$
be the 3-dimensional Lie algebra with Lie brackets 
$[a,b]=0,\ [a,c]=\lambda c,\ [b,c]=\mu c.$

It is easy to see that $[\ad a, \ad b] = 0 $,
and as a consequence $\ad \ln(exp(a)exp(b))
=\ad a+\ad b$. Moreover, for any $\xi(t,u) \in \campo[[t,u]]$ one has
\begin{eqnarray*}
&  \tr(\xi(\ad a, \ad b)) 
= \xi(\lambda,\mu)
+ 2\xi(0,0)
.
\end{eqnarray*}

\end{Remark}

\begin{Theorem}
\label{tm41}
A universal solution  of the  Kashiwara-Vergne conjecture has 
\begin{eqnarray*}
&\rho=0,
\\&
  \gamma(t)
=  \gamma_\alpha(t)
:= \beta_\alpha(t)
- \beta_\alpha(0)
+ \psi'(0)
- \psi'(t)
.
\end{eqnarray*}
\end{Theorem}
\begin{Proof*}
 Let $ \epsilon\in\campo$, 
we have 
$     (\ad Y\epsilon)
\circ \delta_2 B(X,Y\epsilon) 
    = \ad Y
\circ \gamma(\ad X)\epsilon 
    + o(\epsilon)
$ 
and 
$$
 \delta_1	 A(X,Y\epsilon)
= \rho\ id
- \sum_{n\geq 1}
  \beta_n
  \sum_{j=0}^{n-1}
  (\ad X)^j\circ \ad \big((\ad X)^{n-1-j}Y\big)
  \epsilon
+ o(\epsilon)
.
$$
Let $\gG$ be the Lie algebra in Remark \ref{os41}, $X=a$, 
and 
$Y=b$. We get
\begin{eqnarray*}
&&
  \tr( \ad b\circ\delta_2 B(a,b\epsilon) 
     )
  \epsilon
= \tr\big( \ad b\circ \gamma(\ad a)\epsilon + o(\epsilon)
   \big)
,
\\
&& \tr( \ad a\circ\delta_1 A(a,b\epsilon)
    )
= \tr(     \rho\ \ad a
- \sum_{n\geq 1}
  \epsilon\beta_n
  (\ad a)^{n}\circ\ad b
+ o(\epsilon)
),
\\
&&\tr
  \left(
	 \psi(\ad a) 
       + \psi(\epsilon\ \ad b)
       - \psi(\ad a + \epsilon\ \ad  b)
  \right)=
\\
&&\quad \quad \quad \quad \quad \quad 
= \tr
  \left(   -\epsilon(\psi'(\ad a)-\psi'(0)\ad a)\circ\ad b
	 + o(\epsilon)
  \right)
.
\end{eqnarray*}
In particular Equation (\ref{secondaE}) gives
\begin{eqnarray*}
&&
  \tr(\rho\ \ad a)
= 0,  
\\
&& \tr\Big( 
    \big(
     \gamma(\ad a) 
    - \beta(\ad a) +\beta(0)
    + \psi'(\ad a)
    - \psi'(0) 
    \big)
    \circ\ad b
\Big)
= 0
.
\end{eqnarray*}
Rescaling $a$, $b$ and using the properties of $\gG_{\lambda,\mu}$ we get 
\begin{eqnarrayqed*}
&&\rho
= 0, 
\\ 
&&
  \gamma(t) 
- \beta(t)
+ \beta(0)
+ \psi'(t)
- \psi'(0)
= 0
.
\end{eqnarrayqed*}
\end{Proof*}


\begin{Remark}
If one replaces $\psi$ by a series $f= f_0+f_1t+\cdots\in\campo[[t]]$ one 
gets another 
conjecture that 
one can call an $f$-Kashiwara-Vergne conjecture. Then,  

\begin{enumerate}
\item[$i)$]
   Theorem \ref{tm41} is modified by 
replacing $\psi$ by $f$ and  adding $f_0=0$ in the conclusion.
\item[$ii)$] 
  Theorems  \ref{tm31bis}, \ref{tm31}, and part $i)$  imply 
$
  f(t)_{even}
= \psi(t)_{even}
$, 
otherwise  
a universal solution of  the $f$-Kashiwara-Vergne 
conjecture does not exist.
\item[$iii)$]
 If a universal solution of  the $f$-Kashiwara-Vergne 
conjecture has 
$A(X,Y) = B(-Y,-X)$ then one can show that 
$f(t)_{odd} = f_1t$ (we stress that $\psi(t)_{odd}= \psi'(0)t$). 
To get an easy proof one can use the Lie algebra of  Remark \ref{os41}
\end{enumerate}
\end{Remark}


\begin{appendix}
\section{Comparison with  quadratic solutions (by E. Petracci)}

In the previous sections we  did  not determine the  value of the constant 
$\alpha$. Imposing the symmetry condition we obtain
$
  \beta_\alpha(\ad X)Y + o(Y)
= -\alpha Y
- \gamma(-\ad Y)X
+ o(X) 
$, 
so 
$
  \alpha
- \frac{1}{2}
=:  \beta_\alpha(0)~=-\alpha
$. 
Hence  a universal symmetric solution has 
$
  \alpha
= \frac{1}{4}
$.

\vskip0.5cm

Vergne and Alekseev-Meinrenken both considered a 
quadratic Lie algebra  and obtained  symmetric solutions.  
It is natural to ask whether these solutions are universal. In fact,  
 quadratic Lie algebras have the special property 
$\tr((\ad X)^{2n}\circ\ad Y) = 0$ for any $n\in\Naturali$ 
and any couple of vectors $X,Y$,   
which 
  simplifies the equation with traces 
(\ref{secondaE}). 

We have seen that  a 
universal symmetric solution of the Kashiwara-Vergne conjecture 
has
\begin{eqnarray*}
&&B(X,Y) 
= \frac{1}{4} 
  X 
+ \Big(
   \underbrace{
	  \beta_{\frac{1}{4} }(\ad X)
	- \psi'(\ad X)
	+ \frac{1}{2}id 
   }_{\gamma_\frac{1}{4}(\ad X)}
  \Big)
  (Y) 
+ o(Y)
,
\\
&&A(X,Y) 
= \beta_{\frac{1}{4} }(\ad X)(Y )
+ \frac{1}{2}
  (\pi_{\gamma_\frac{1}{4}(t)}(t,u):[Y,Y])_X
+ o(Y^2)
\end{eqnarray*}
with $\beta_\alpha(t)$ given in Theorem \ref{tm31bis}, and 
$\pi_{\gamma_\frac{1}{4}(t)}(t,u)$ given in Theorem \ref{tm31}.

\begin{Remark} {\bf\small (Vergne's solution for quadratic Lie algebras)}
\label{rem51}
\newline
We denote by $B_V(X,Y)$ the $B$ found by M. Vergne in  {\bf \cite{Ver}}. 
Following her 
paper we 
find  
$
  B_V(X,Y)
= \frac{1}{4}X 
+ \gamma_V(\ad X)(Y)
+ o(Y)
$. 
Let  
$$
   R(t) 
:= \frac{e^t-e^{-t}-2t}{t^2}
,
$$
 after a bit long calculation 
we see that the series $\gamma_V(t)$ is 
given 
by 
\begin{eqnarray*}
  t{\gamma_V}'(t)
+ 2{\gamma_V}(t)
&=& \frac{1}{8}t
- \frac{1}{2}
  t\varphi_1(-t)R(t){\varphi_1}'(t)=
\\
&=& \frac{1}{8}t
  + \frac{1}{12}t^2
  + \frac{1}{72}t^3 
  - \frac{1}{360}t^4 
+o(t^4)
.
\end{eqnarray*}
This differential equation  gives $\gamma_V(t)_{odd} = 
\gamma_\frac{1}{4}(t)_{odd}$. 
A universal solution has 
$$
  t\gamma_\frac{1}{4}'(t)
+ 2\gamma_\frac{1}{4}(t)
= \frac{1}{8}t 
+ \frac{1}{12}t^2
+ \frac{1}{72}t^3 
- \frac{1}{480}t^4 
+ o(t^4)
.
$$ 
In particular the symmetric solution found by M. Vergne for a quadratic 
Lie algebra 
is not universal. 
\end{Remark}

\begin{Remark} 
\label{rem52}
{\bf\small(Alekseev-Meinrenken's solution for quadratic Lie 
algebras)}
\newline
Let 
$  B_{AM}(X,Y) 
= \alpha_{AM}
+ \gamma_{AM}(\ad X)(Y)
+ o(Y) 
$  
the $B$ found by Alekseev and Meinrenken, $\beta_{AM}(t)$ their series 
$\beta(t)$, etc.

  Following 
their paper {\bf\cite{AlM}} and the paper {\bf\cite{Ver}} of  M. Vergne,   
after some efforts we find  the following formulas. Let 
$g(t)=\frac{1}{2}R(t)$, and 
$\Pi(t)$ be the formal power series 
such that \ \ \ 
$ t\Pi'(t)
+ 2\Pi(t) 
= \frac{1}{2}
  \varphi_1(t)^{-1}
- g(t)
  \varphi_1(-t)
  (1-\varphi_1(t))  
.
$  
Then 
\begin{eqnarray*}
 && \beta_{AM}(t)
= \Pi(t)
- \frac{1}{4}
  \left( g(t)\varphi_1(-t)-\frac{1}{2}\varphi_1(t)^{-1} 
  \right)
  t
- \frac{1}{2}
  \varphi(t)^{-1}+ 
\\
&&
\quad\quad\quad
\quad\quad\quad
\quad\quad\quad
+ g(t)
  \varphi_1(-t)
  (1-\varphi_1(t))
,
\\
&&  \rho_{AM}
= 0
,
\\   
&&  \gamma_{AM}(-t)+\gamma_V(t)
= \varphi_1(-t)
  g(t)
  t
  \left( \beta_{AM}(-t)-\frac{1}{4}\varphi_1(t)-\varphi_1'(t)
\right)+
\\	
&&
\quad\quad\quad
\quad\quad\quad
\quad\quad\quad
- \frac{1}{2}
  \varphi_1(t)^{-1}
  t
  \beta_{AM}(-t)
- \frac{1}{8}t
,
\\
&&  \alpha_{AM}
= \frac{1}{4}
.
\end{eqnarray*}
Using Maple we get 
$\beta_{AM}(t)= \beta_{\frac{1}{4}}(t)$,  
$\gamma_{AM}(t)_{odd}=\gamma_V(t)_{odd}$, and 
\begin{eqnarray*}
  t\gamma_{AM}'(t)
+ 2
  \gamma_{AM}(t)
= \frac{1}{8}
  t
+ \frac{1}{12}
  t^2
+\frac{1}{72}
  t^3
- \frac{1}{720}
  t^4
+ o(t^4)
.
\end{eqnarray*}
In particular  the symmetric solution of Alekseev and Meinrenken is not 
universal, and it is 
different from 
the solution of Vergne.
\end{Remark}
\noindent

\end{appendix}


\references
\lastpage
\end{document}